\documentclass[11pt]{amsart}
\usepackage[colorlinks=true, pdfstartview=FitV, linkcolor=blue, 
citecolor=blue]{hyperref}

\usepackage{amssymb,amsmath,amsfonts}
\usepackage{bm}
\usepackage{enumerate}
\usepackage{array}

\allowdisplaybreaks

\newcommand{\NN}{\mathbb{N}}
\newcommand{\RR}{\mathbb{R}}

\newcommand{\ZZ}{\mathbb{Z}}
\newcommand{\QQ}{\mathbb{Q}}
\newcommand{\rar}{\rightarrow}
\newcommand{\mc}{\mathcal}

\newcommand{\bx}{\mathbf{x}}

\newtheorem*{theorem*}{Theorem}
\newtheorem{theorem}{Theorem}

\begin{document}

\title[Badly approximable points in solenoids]{Badly approximable points for\\diagonal approximation in solenoids}
\author{Huayang Chen, Alan Haynes}

\thanks{AH: Research supported by NSF grant DMS 2001248.\\
\phantom{A..}MSC 2020: 11J61, 11J83}

\keywords{Badly approximable numbers, Diophantine approximation, $p$-adic numbers}

\begin{abstract}
In this paper we investigate the problem of how well points in finite dimensional $p$-adic solenoids can be approximated by rationals. The setting we work in was previously studied by Palmer, who proved analogues of Dirichlet's theorem and the Duffin-Schaeffer theorem. We prove a complementary result, showing that the set of badly approximable points has maximum Hausdorff dimension. Our proof is a simple application of the elegant machinery of Schmidt's game.
\end{abstract}

\maketitle

\section{Introduction}
There are several natural ways of extending the classical study of Diophantine approximation to spaces of $p$-adic numbers. The first and most obvious, which was pioneered by K.~Mahler \cite{Mahl1940}, Jarn\'{i}k \cite{Jarn1945}, and Lutz \cite{Lutz1955}, is to study approximations in $\QQ_p$ in the usual sense by simply restricting to $\ZZ_p$ (see also \cite{BereBernKova2005,BernDickYuan1999,HaynMund2015,Hayn2010,Li2016,Liu2019,Swee1974}). A second way, which was introduced by Choi and Vaaler \cite{ChoiVaal1999} and studied using techniques from the geometry of numbers over the adeles \cite{BombVaal1983,BurgVaal1993}, is to work in projective space over $\QQ_p$ (see also \cite{GhosHayn2016,HarrHuss2017}). A third way, studied in some form by D.~G.~Cantor \cite{Cant1965}, and more recently by Palmer \cite{Palm2015}, is to work in a $p$-adic solenoid, which is the quotient of $\RR\times\QQ_p$ by a diagonally embedded lattice (see also \cite{FurnHaynKoiv2019,HaynMund2013}).

Analogues of Dirichlet's theorem, Khintchine's theorem \cite{Khin1924,Khin1926}, and the Duffin-Schaeffer theorem \cite{DuffScha1941} have been proven in all three of the above mentioned settings (see \cite{Hayn2010}, \cite{GhosHayn2016}, and \cite{Palm2015}, respectively). However, complementary results for sets of badly approximable $p$-adic numbers (i.e. analogues of Jarn\'{i}k's classical results \cite{Jarn1928} for the reals) have only been studied in the first two settings: in $\ZZ_p$ it was shown by Abercrombie \cite{Aber1995} that the set of badly approximable numbers has full Hausdorff dimension (see also \cite[Section 5.4]{KrisThorVela2006}), while for projective space over $\QQ_p$ the analogous result was recently established by Harrap and Hussain \cite{HarrHuss2017}. Our goal in this article is to complete this part of the story by determining the size of the set of badly approximable numbers in the third setting of diagonal approximation.

Let $\mc{P}=\{p_1,\ldots ,p_k\}$ be a finite set of distinct prime numbers and define
\begin{equation*}
\QQ_\mc{P}=\RR\times\QQ_{p_1}\times\cdots\times\QQ_{p_k}.
\end{equation*}
We consider $\QQ_\mc{P}$ as a topological group with the usual product topology inherited from its factors. Writing $\bx\in\QQ_\mc{P}$ as $\bx=(x_\infty,x_{p_1},\ldots ,x_{p_k})$, this topology is the metric topology arising from the sup-norm
\begin{equation*}
|\bx|=\max\left\{|x_\infty|_\infty,|x_{p_1}|_{p_1},\ldots ,|x_{p_k}|_{p_k}\right\}.
\end{equation*}
The diagonal embedding
\begin{align*}
\iota:\QQ\rar\QQ_\mc{P}
\end{align*}
defined by
\begin{align*}
\gamma\mapsto (\gamma,\ldots,\gamma)
\end{align*}
is a homomorphism of additive groups. For $\gamma\in\QQ$ we will write $\bm{\gamma}=\iota(\gamma)$. The group
\begin{equation*}
\Gamma_\mc{P}=\iota\left(\ZZ\left[\frac{1}{p_1\cdots p_k}\right]\right)
\end{equation*}
is a discrete and co-compact subgroup of $\QQ_\mc{P}$. It is not difficult to show (see \cite[p.111]{Palm2015}) that the set
\begin{equation*}
\mc{F}=[0,1)\times\ZZ_{p_1}\times\cdots \times\ZZ_{p_k}
\end{equation*}
is a strict fundamental domain for the quotient $\QQ_\mc{P}/\Gamma_\mc{P}$. Our problem of interest, which is derived by direct analogy with the classical theory of approximation in $\RR/\ZZ$, is to determine how well elements of $\mc{F}$ can be approximated by diagonally embedded fractions of the form $\bm{\beta}/\gamma$, where $\beta,\gamma\in\ZZ[1/(p_1\cdots p_k)]$. A first result in this direction is the following version of Dirichlet's theorem, proved by Palmer.
\begin{theorem}\cite[Theorem 1.1]{Palm2015}\label{thm.DTSolenoid}
Let $M=\max\{p_1,\ldots,p_k\}$. For any $\bm{x}\in\QQ_\mc{P}$ and for any $N\in\NN$, there exist $\beta,\gamma\in\ZZ[1/(p_1\cdots p_k)]$ satisfying $0<|\bm{\gamma}|\le N$ and
\begin{equation*}
|\gamma\bm{x}-\bm{\beta}|\le\frac{M}{N}.
\end{equation*}
\end{theorem}
It also follows from \cite[Theorem 1.2]{Palm2015} that, for any $\epsilon>0$, the set of $\bm{x}\in\mc{F}$ for which there are infinitely many pairs $\beta,\gamma\in\ZZ[1/(p_1\cdots p_k)], \gamma\not=0,$ satisfying
\begin{equation}\label{eqn.ApproxIneq}
|\gamma\bm{x}-\bm{\beta}|\le\frac{\epsilon}{|\bm{\gamma}|},
\end{equation}
is a set of Haar measure $0$. We therefore define the set of \textit{badly approximable numbers in $\QQ_\mc{P}/\Gamma_\mc{P}$} to be the set of $\bm{x}\in\mc{F}$ for which there exists a constant $\delta=\delta(\bm{x})>0$ such that
\begin{equation*}
|\gamma\bm{x}-\bm{\beta}|\ge\frac{\delta}{|\bm{\gamma}|},
\end{equation*}
for all $\beta,\gamma\in\ZZ[1/(p_1\cdots p_k)]$ with $\gamma\not=0$. We denote this set by $\mc{B}_\mc{P}$. By what we have already said, it follows that $\mc{B}_\mc{P}$ is a set of Haar measure $0$. However, our main result in this paper is that from the point of view of Hausdorff dimension, this set is still large.
\begin{theorem}\label{thm.HDBad}
With respect to the sup-norm metric on $\QQ_\mc{P}$, we have that
\begin{equation*}
\dim_{H}\mc{B}_\mc{P}=\dim_H \QQ_\mc{P}=k+1.
\end{equation*}
\end{theorem}
To prove this theorem directly using the definition and basic theory of Hausdorff dimension does not appear to be very simple. One difficulty is that sets defined by inequalities of the form \eqref{eqn.ApproxIneq} (which must be removed during the construction of $\mc{B}_\mc{P}$) are not balls in the sup-norm metric. Instead, they are the sets defined by the inequalities \eqref{eqn.RegAvoid} below, and their diameters in each of the Archimedean and $p$-adic directions depend individually on the Archimedean and $p$-adic absolute values of $\gamma$, which may be different. However, as we will see, viewing the problem in the framework of Schmidt's game allows us to overcome this difficulty.

We will briefly review the definitions and relevant results about Schmidt's game in Section \ref{sec.SchmGame}, and we will give the proof of Theorem \ref{thm.HDBad} in Section \ref{sec.PfOfHDBad}.

\section{Schmidt's game}\label{sec.SchmGame}
The game we are going to describe was first introduced by Wolfgang Schmidt in \cite{Schm1966}, and it is now referred to as Schmidt's game. Let $X$ be a complete metric space with metric $\mathrm{d}$, and for $x\in X$ and $\rho\ge 0$ let $B(x,\rho)$ denote the closed ball of radius $\rho$ centered at $x$. Define a partial order $\prec$ on the set $X\times\RR^+$ by the rule
\begin{equation*}
(x_1,\rho_1)\prec (x_2,\rho_2) \quad\text{if and only if}\quad \rho_1+\mathrm{d}(x_1,x_2)\le\rho_2.
\end{equation*}
If $(x_1,\rho_1)\prec (x_2,\rho_2)$ then it follows from the triangle inequality that 
\begin{equation}\label{eqn.BallCont}
B(x_1,\rho_1)\subseteq B(x_2,\rho_2).
\end{equation}
However, it is important in what follows to understand that the converse of this statement is not true in general (although it is true when $X$ is Euclidean space). In particular, when $X=\QQ_\mc{P}$, it can happen that \eqref{eqn.BallCont} holds but that $(x_1,\rho_1)\not\prec (x_2,\rho_2)$.

Now let $\alpha$ and $\beta$ be real numbers satisfying $0<\alpha,\beta<1$, and let $S\subseteq X.$ Two players, Bob and Alice, alternate choosing a sequence of balls in $X$. Bob  plays first and chooses any ball $B_0=B(b_0,\rho_0)$. Then Alice chooses a ball $A_0=B(a_0,\rho_0')$ satisfying $\rho_0'=\alpha\rho_0$ and $(a_0,\rho_0')\prec (b_0,\rho_0)$. Next, Bob chooses a ball $B_1=B(b_1,\rho_1)$ satisfying $\rho_1=\beta\rho_0'$ and $(b_1,\rho_1)\prec (a_0,\rho')$. Then Alice chooses a ball $A_1=B(a_1,\rho_1')$ satisfying $\rho_1'=\alpha\rho_1$ and $(a_1,\rho_1')\prec (b_1,\rho_1)$, and so on, creating a sequence of nested balls
\begin{equation*}
B_0\supseteq A_0\supseteq\cdots\supseteq B_n\supseteq A_n\supseteq\cdots,
\end{equation*}
with radii tending to $0$ and satisfying $\rho_n'=\alpha\rho_n$ and $\rho_{n+1}=\beta\rho_n'$ for each $n\in\NN$. Since this game is played on a complete metric space, the intersection of these balls is a single point $x_\infty\in X$. We say that $S$ is an \textit{$(\alpha,\beta)$-winning set} if Alice can always choose her balls so that $x_\infty\in S$, regardless of how Bob plays. We say that $S$ is an \textit{$\alpha$-winning set} if it is $(\alpha,\beta)$-winning for all $\beta\in(0,1)$, and we say that $S$ is \textit{winning} if it is $\alpha$-winning for some $\alpha$.

Since Bob has free choice of $B_0$, it is obvious that any $(\alpha,\beta)$-winning set must be dense in $X$. In fact, Schmidt proved much stronger conclusions about such sets. For example, he showed that any countable intersection of $\alpha$-winning sets in $X$ is also $\alpha$-winning \cite[Theorem 2]{Schm1966}, and he also gave non-trivial lower bounds for Hausdorff dimensions of $(\alpha,\beta)$-winning sets \cite[Theorem 6]{Schm1966}, in the case when $X$ is a Hilbert space. For our application, we require an analogue of the latter result for the case when $X=\QQ_\mc{P}$.

For completeness, recall that for a metric space $(X,\mathrm{d})$, a set $S\subseteq X$, and a real number $s\ge 0$, the $s$-dimensional Hausdorff outer measure of $S$ is defined as
\begin{equation*}
\mathcal{H}^s(S)=\lim_{\rho\rar 0^+}\left(\inf\left\{\sum_{i=1}^\infty\rho_i^s:\{B(x_i,\rho_i)\}_{i\in\NN} \text{ is a } \rho\text{-cover of } S\right\}\right),
\end{equation*}
where a $\rho$-cover $\{B(x_i,\rho_i)\}_{i\in\NN}$ of $S$ is any cover of $S$ by balls satisfying $0<\rho_i<\rho$ for all $i$. The Hausdorff dimension of $S$ is then defined as
\begin{equation*}
	\dim_H S=\inf\{s\ge 0:\mathcal{H}^s(S)=0\}.
\end{equation*}

Let $X=\QQ_\mc{P}$ and, for $\beta\in (0,1)$, let $N(\beta)$ be the largest integer with the property that, for every $\rho>0$, every ball $B'=B(\bm{x},\rho)$ in $X$ contains a collection of $N=N(\beta)$ balls $B_1',\ldots ,B_N'$, with centers $\bm{x}^{(i)}$ and radii $\beta\rho$, with pairwise disjoint interiors, and satisfying
\begin{equation}\label{eqn.PrecCond}
(\bm{x}^{(i)},\beta\rho)\prec (\bm{x},\rho)
\end{equation}
for each $1\le i\le N$. Then, by a slight modification of the proof of \cite[Theorem 4.1]{Kris2006}, it can be shown that the Hausdorff dimension of any $(\alpha,\beta)$-winning set $S\subseteq X$ must satisfy the inequality
\begin{equation}\label{eqn.HDBd}
\dim_H S\ge\frac{\log N(\beta)}{|\log\alpha\beta|}.
\end{equation}
To establish \eqref{eqn.HDBd}, the only difference from the proof of \cite[Theorem 4.1]{Kris2006} is that we do not have the ``ball intersection property'' in $\QQ_\mc{P}$. However, in place of this it is sufficient to observe that if $r_1<r_2$ then any ball of radius $r_1$ in $\QQ_\mc{P}$ cannot intersect more than two elements from any collection of disjoint balls of radius $r_2$ in $\QQ_\mc{P}$. 

Finally, to bound $N(\beta)$, write $\mc{P}=\{p_1,\ldots ,p_k\}$ and observe that \eqref{eqn.PrecCond} is satisfied if and only if
\begin{equation*}
\max\left\{|x^{(i)}_\infty-x_\infty|_\infty,|x^{(i)}_{p_1}-x_{p_1}|_{p_1},\ldots ,|x^{(i)}_{p_k}-x_{p_k}|_{p_k}\right\}\le\rho (1-\beta).
\end{equation*}
This gives $\lfloor 1/\beta\rfloor$ choices for the Archimedean component, each of which will give rise to non-overlapping Archimedean intervals of radius $\beta\rho$. For each prime $p_j$, as long as $1-\beta>1/p_j$, the number of choices for $p_j$-adic components (which will satisfy the above inequality and also give rise to non-overlapping $p_j$-adic balls of radius $\beta\rho$) is at least
\begin{equation*}
p_j^{-1}p_j^{\lfloor\log_{p_j}(1/\beta)\rfloor}.
\end{equation*}
Combining these estimates, as long as $0<\beta<1/2$, we have that
\begin{equation}\label{eqn.N(beta)Bd}
N(\beta)\ge\lfloor1/\beta\rfloor\prod_{j=1}^kp_j^{-1}p_j^{\lfloor\log_{p_j}(1/\beta)\rfloor}\ge\frac{1}{2(p_1\cdots p_k)^2\beta^{k+1}}.
\end{equation}
If $S\subseteq\QQ_\mc{P}$ is $\alpha$-winning then applying the bounds \eqref{eqn.HDBd} and \eqref{eqn.N(beta)Bd} and taking the limit as $\beta\rar 0^+$, we obtain
\begin{equation*}
\dim_H S=k+1.
\end{equation*}
It is also not difficult to show that the Hausdorff dimension of $\QQ_\mc{P},$ with respect to the metric we are using, is $k+1$. Therefore, to prove Theorem \ref{thm.HDBad} it is sufficient to show that $\mc{B}_\mc{P}$ is an $\alpha$-winning set, for some $0<\alpha<1$.

\section{Proof of Theorem \ref{thm.HDBad}}\label{sec.PfOfHDBad}

Our proof is quite similar to the proof of \cite[Theorem 3]{Schm1966}, albeit with differences to account both for the non-Archimedean directions and for the fact that the sets we need to avoid are not balls in the sup-norm metric on $\QQ_\mc{P}$. Take
\begin{equation*}
\alpha=\min_{1\le i\le k}\frac{1}{p_i^{2}},
\end{equation*}
let $\beta\in (0,1)$, and set $c=c(\alpha,\beta)=1+\alpha\beta-2\alpha$. We will show that Alice can always play in a way so that the intersection point $\bm{x}_\infty\in\QQ_\mc{P}$ satisfies
\begin{equation*}
|\gamma\bm{x}_\infty-\bm{\beta}|\ge\frac{\delta}{|\bm{\gamma}|}
\end{equation*}
for all $\beta,\gamma\in\ZZ[1/(p_1\cdots p_k)]$, $\gamma\not=0,$ with
\begin{equation*}
\delta=\min\left\{\alpha,c/2\right\}\cdot\min\left\{\rho_0,\alpha^2\beta^2c/8\right\}.
\end{equation*}
Note that since $0<\alpha<1/2,$ we always have that $0<c<1$.

First of all, if $\rho_0> \alpha\beta c/8$ then there will be an integer $n\ge 1$ for which
\begin{equation*}
\alpha^2\beta^2 c/8<\rho_n\le \alpha\beta c/8.
\end{equation*}
In this case we will have that $\rho_0>\rho_n>\alpha^2\beta^2c/8$, and Alice can play as if $B_n$ were Bob's first ball, without affecting the choice of $\delta$ in the argument. Therefore, we will assume without loss of generality that $\rho_0\le\alpha\beta c/8$.

Let $t\ge 1$ be the unique integer with the property that 
\begin{equation*}
\alpha\beta c/2\le (\alpha\beta)^t<c/2
\end{equation*}
and let $R>0$ be defined by
\begin{equation*}
\frac{1}{R^2}=(\alpha\beta)^{t}.
\end{equation*}
We will prove by induction that, for every integer $n\ge 0$, Alice can play in a way that for every $\bm{x}\in B_{nt+1}$ and for every  $\beta,\gamma\in\ZZ[1/(p_1\cdots p_k)]$ with $0<|\bm{\gamma}|< R^n$, we have that
\begin{equation}\label{eqn.BallStrat1}
|\gamma\bm{x}-\bm{\beta}|\ge\frac{\delta}{|\bm{\gamma}|}.
\end{equation}
When $n=0$ this is trivially true, since $|\bm{\gamma}|\ge 1$ for all non-zero $\gamma\in\ZZ[1/(p_1\cdots p_k)]$. Suppose it holds for all non-negative integers less than $n$, for some $n\in\NN$. Then Alice needs to choose the balls $A_{(n-1)t+\ell}$, with $1\le \ell\le t$, so that no matter how Bob plays, equation \eqref{eqn.BallStrat1} is guaranteed to hold for all $\bm{x}\in B_{nt+1}$ and for all $\beta$ and $\gamma$ with $R^{n-1}\le|\bm{\gamma}|<R^n.$

A key observation is that, for each $n$, there is at most one possible pair $\beta,\gamma$ which Alice needs to avoid. To see why this is true, suppose that $\beta,\gamma,\beta',\gamma'\in\ZZ[1/(p_1\cdots p_k)]$ and $\bm{x},\bm{x}'\in B_{nt+1}$ are chosen so that $R^{n-1}\le |\bm{\gamma}|,|\bm{\gamma}'|<R^n$ and
\begin{equation*}
|\gamma\bm{x}-\bm{\beta}|<\frac{\delta}{|\bm{\gamma}|}\qquad\text{and}\qquad|\gamma'\bm{x}'-\bm{\beta}'|<\frac{\delta}{|\bm{\gamma}'|}.
\end{equation*}
Then we have that
\begin{align*}
|\gamma'\bm{\beta}-\gamma\bm{\beta}'|&=|\gamma'\bm{\beta}-\gamma'\gamma\bm{x}+\gamma'\gamma\bm{x}-\gamma\gamma'\bm{x}'+\gamma\gamma'\bm{x}'-\gamma\bm{\beta}'|\\
&<\frac{\delta|\bm{\gamma}'|}{|\bm{\gamma}|}+2|\bm{\gamma}'||\bm{\gamma}|\rho_{(k-1)t+1}+\frac{\delta|\bm{\gamma}|}{|\bm{\gamma}'|}\\
&\le 2\delta R+2\rho_0 R^2<4\rho_0R^2\le \frac{\alpha\beta c(\alpha\beta)^{-t}}{2}\le 1.
\end{align*}
If $\beta/\gamma\not=\beta'/\gamma'$ then, since 
\begin{equation*}
|\gamma'\beta-\gamma\beta'|_\infty\prod_{i=1}^k|\gamma'\beta-\gamma\beta'|_{p_i}\in\ZZ,
\end{equation*}
we would have that
\begin{equation*}
|\gamma'\bm{\beta}-\gamma\bm{\beta}'|\ge 1.
\end{equation*}
This would give a contradiction, therefore we conclude that $\beta/\gamma=\beta'/\gamma'$.

Now suppose that $\beta,\gamma$ satisfy $R^{n-1}\le |\bm{\gamma}|<R^n$ and that there is a point $\bm{x}\in B_{(n-1)t+1}$ with
\begin{equation*}
|\gamma\bm{x}-\bm{\beta}|<\frac{\delta}{|\bm{\gamma}|}.
\end{equation*}
Note that the set of all such points is precisely the set of $\bm{x}\in B_{(n-1)t+1}$ which satisfy
\begin{align}\label{eqn.RegAvoid}
\left|x_\infty-\frac{\beta}{\gamma}\right|_\infty<\frac{\delta}{|\gamma|_\infty|\bm{\gamma}|}\qquad\text{and}\qquad \left|x_{p_i}-\frac{\beta}{\gamma}\right|_{p_i}<\frac{\delta}{|\gamma|_{p_i}|\bm{\gamma}|},
\end{align}
for each $1\le i\le k$. Thus, in order to avoid the problematic region, we need to ensure that one of the above inequalities is violated. We distinguish two cases, depending on whether or not the Archimedean absolute value of $\gamma$ is largest.

\noindent Case 1: If $|\bm{\gamma}|=|\gamma|_\infty$ then in the Archimedean direction Alice needs to avoid an interval of radius at most $\delta/R^{2(n-1)}.$ If the center of this interval is in the left half of the Archimedean component of $B_{(n-1)t+1}$ then Alice should always choose her balls so that their Archimedean components are as far to the right as possible, otherwise she should always choose them so that they are as far to the left as possible. At each step, she should choose all of the $p_i$-adic components of the centers of her balls to be the same as Bob's, thus guaranteeing that
\begin{equation*}
(\bm{a}_{(n-1)t+\ell},\rho_{(n-1)t+\ell}')\prec (\bm{b}_{(n-1)t+\ell},\rho_{(n-1)t+\ell})
\end{equation*}
for each $1\le\ell\le t$. It then follows from the argument in Schmidt's paper \cite[Section 7]{Schm1966} that Alice has a winning strategy. We note for completeness that Schmidt's argument requires that $c>0$ and also that 
\begin{equation*}
\delta\le (c/2)\min\left\{\rho_0,\alpha^2\beta^2c/8\right\},
\end{equation*}
which is part of the reason for our choice of $\delta$.

\noindent Case 2: If $|\bm{\gamma}|=|\gamma|_{p_i}$ for some $1\le i\le k$ then we have that
\begin{equation*}
\frac{\delta}{|\gamma|_{p_i}|\bm{\gamma}|}\le \frac{\delta}{R^{2(n-1)}}\le \alpha\rho_{(n-1)t+1}\le\frac{\rho_{(n-1)t+1}}{p_i^2}.
\end{equation*}
The ball in $\QQ_{p_i}$ of radius $\rho_{(n-1)t+1}/p_i$, centered at the $p_i$-adic component of $\bm{b}_{(n-1)t+1}$, is a disjoint union of $p_i$ balls of radius $\rho_{(n-1)t+1}/p_i^2$, and only one of these can intersect the ball of radius $\delta/(|\gamma|_{p_i}|\bm{\gamma}|)$ centered at $\beta/\gamma$. This leaves Alice with $p_i-1$ choices for the $p_i$-adic component of $\bm{a}_{(n-1)t+1}$, all of which avoid the point $\beta/\gamma$. Assuming she chooses the rest of the components of $\bm{a}_{(n-1)t+1}$ to exactly match Bob's, we have that
\begin{equation*}
\rho_{(n-1)t+1}'+|\bm{a}_{(n-1)t+1}-\bm{b}_{(n-1)t+1}|\le\frac{\rho_{(n-1)t+1}}{p_i^2}+\frac{\rho_{(n-1)t+1}}{p_i}<\rho_{(n-1)t+1},
\end{equation*}
which guarantees that
\begin{equation*}
(\bm{a}_{(n-1)t+1},\rho_{(n-1)t+1}')\prec (\bm{b}_{(n-1)t+1},\rho_{(n-1)t+1}).
\end{equation*}
From that point on, for the rest of the balls $A_{(n-1)t+\ell}$ with $2\le\ell\le t$, Alice is free to play however she wants. This is because once the ball around the point $\beta/\gamma$ has been avoided once in the $p_i$-adic component, it cannot intersect any of the future balls. This completes the proof of the theorem.

\vspace{.15in}

{\footnotesize
\noindent
Department of Mathematics,\\
University of Houston,\\
Houston, TX, United States.\\
hchen60@uh.edu, haynes@math.uh.edu

}

\end{document}